\documentstyle{amsppt}
\input epsf.sty
\NoRunningHeads
\TagsOnRight

\def\l{\lambda}
\def\b{\beta}

\def\mathbb{\Bbb}
\def\C{\Bbb C}
\def\mathfrak{\goth}
\def\w{^\wedge}
\def\la{\lambda}

\def\net{\varnothing}
\def\Hom{\text{\rm Hom}}
\def\End{\text{\rm End}\,}
\def\ZMN{Z(M,N)}
\def\CSn{\Bbb C[S_n]}
\def\I{\Bbb I}
\def\al{\alpha}

\def\la{\lambda}
\def\ar{\nearrow}
\def\T{\operatorname{Tab}}
\def\Y{\mathbb Y}
\def\ssim{\approx}
\def\partn{\vdash}
\def\Q{\Bbb Q}
\def\R{\Bbb R}
\def\Z{\Bbb Z}

\topmatter

\title{A NEW APPROACH TO THE REPRESENTATION THEORY OF
THE SYMMETRIC GROUPS. II}
\endtitle
\author{A.~M.~Vershik and A.~Yu.~Okounkov}%{512.547}
\endauthor

\abstract{The present paper is a revised Russian translation of the paper
``A new approach to representation
theory of symmetric groups,'' {\it Selecta Math., New Series}, {\bf 2}, No.~{\rm 4,
581--605 (1996)}. Numerous modifications to the text were made by the first
author for this publication. Bibliography: $35$ titles.}
\endabstract
\endtopmatter

\document

\hfill{\bf To the memory of D.~Coxeter}

\head Preface
\endhead

This paper is a revised Russian translation of a paper by the same
authors (see the reference below) and is devoted to a nontraditional
approach to the representation theory of the symmetric groups (and, more
generally, to the representation theory of Coxeter and local groups).
The translation was prepared for the Russian edition of the book
W.~Fulton, {\it Young Tableaux. With Applications to Representation
Theory and Geometry}, Cambridge Univ.\ Press, Cambridge, 1997,
which, hopefully, will appear sooner or later. In the editor's preface to the
Russian translation of the book it is explained what is the drawback of
the conventional approach to the representation theory of the
symmetric groups: it does not take into account important properties
of these groups, namely, that they are Coxeter groups, and that they form an
inductive chain, which implies that the theory must be constructed
inductively. A direct consequence of these drawbacks is, in particular,
that Young diagrams and tableaux appear {\it ad hoc}; there
presence in the theory is justified only after the proof of the branching theorem.

The theory described in this paper is intended to correct these defects.
The first attempt in this direction was the paper [30]
by the first author,
 in which it was proved that if we assume  that the
branching graph of irreducible complex representations of the
symmetric groups  is distributive, then it must be the Young graph. As it turned out, this
{\it a priori} assumption
is superfluous --- the distributivity follows directly from the fact
that $S_n$ is a Coxeter group if we involve remarkable generators
of the Gelfand--Tsetlin subalgebra of the group algebra $\Bbb C[S_n]$,
namely, the
Young\footnote{The so-called Young's orthogonal and seminormal form
for the action of the Coxeter transpositions in irreducible
representations were defined in Young's last papers; apparently,
 he regarded
them only as an illustration;
these forms play an essential role in our theory
(see Secs.\ 3, 7). Some time ago A.~Lascoux observed that these generators
were mentioned explicitly in  Young's paper,
and recently
R.~Stanley gave a precise reference. But apparently Young himself
underestimated their importance.}--Jucys\footnote{A.-A.~A.~Jucys (1936--1998)
is a Lithuanian mathematician. The paper \cite{19}, where he introduced
these generators, remained unnoticed for a long time; an English
mathematician G.~Murphy rediscovered them and then found Jucys'
paper.}--Murphy generators (see \cite{19, 30}).
But all numerous later expositions, including the very good
book by Fulton, followed the classical version of the theory, which goes back to
Frobenius, Schur, and Young; although  some nice
simplifications were made, such as von Neumann's lemma, Weyl's lemma,
the notion of tabloids, etc., but the general scheme
of the construction of the theory remained the same.\footnote{Our approach to the representation
theory of the symmetric groups was recently used in \cite{35}.} The reader can find references to the books
on the representation theory of the symmetric groups in the
monograph by James and Kerber \cite{18}, in the book by James \cite{17},
which was
translated into Russian, and in earlier textbooks.

The key point of our approach, which explains
the appearance of Young tableaux as well
as the general idea of our method, is that the points of the spectrum
of the Gelfand--Tsetlin algebra with respect to the Young--Jucys--Murphy
generators are so-called content vectors, i.e.,
 integer vectors  in $\Bbb R^n$
that satisfy certain simple conditions, which follow from the Coxeter relations,
and the coordinates of these integer vectors are the so-called
contents of the boxes of Young tableaux (see Sec.~6);  since the content
vector uniquely determines a Young tableau, it follows that {\it the points
of the spectrum are precisely Young tableaux}. The corresponding eigenvectors
determine  a basis in each representation, and the set of vectors
corresponding to tableaux with a given diagram form a basis of the
irreducible representation of $S_n$ (the Young--Gelfand--Tsetlin basis).
{\it Thus the correspondence ``diagrams'' $\leftrightarrow$ ``irreducible
representations'' obtains a natural (one might say, spectral)
explanation.}

Our approach not only helps to improve
the exposition of classical
results, it also allows us to consider representations of more general
groups and algebras, for example, ``local groups and algebras'' in the
sense of \cite{30}, provided that the group is finite or
the algebra is finite-dimensional. An attempt to apply this method to
other groups and, in particular, to the Coxeter groups of series B--C--D, is
contained in \cite{12} and \cite{28}.

Recently died a distinguished and original mathematician
Donald Coxeter (1907--2003), to whom modern mathematics
owes important and deep ideas and very beautiful
geometric and group constructions. This revised version of
the paper is dedicated to the memory of D.~Coxeter.

\bigskip

\hfill A.~Vershik

\head
0. Introduction
\endhead

The aim of this
paper is to present a new, simple and direct,
approach to the representation theory of
the permutation group $S_n$.

Basically, there are two ways to construct
irreducible complex
representations of $S_n$. The first one is essentially based
on the representation theory of the full linear group $GL(N)$
and the duality between $S_n$ and $GL(N)$ in
the space
$$
\underbrace{\C^N\otimes\C^N\otimes\dots\otimes\C^N}
_{\text{$n$ times}} \,,
$$
which is called the Schur--Weyl duality (see \cite{1}).
The Schur functions, which are characters of
$GL(N)$, play the key role in this
approach. A description of the characters of $S_n$ that is based on the
Schur functions and is close to the original construction by Frobenius
can be found, for example, in \cite{23}.

The other way, usually attributed
to Young with later contributions by von Neumann and Weyl, is based on the
combinatorics of tableaux. In this approach, an irreducible
representation (sometimes called a Specht module) arises as
the unique common component of two simple
representations induced from one-dimensional representations
(the identity representation and the sign
representation) of the same Young subgroup. It is this irreducible component
that one associates with the partition (diagram) corresponding to the
Young subgroup. Since the decomposition of induced representations into
irreducible ones is rather complicated and nonconstructive, the
correspondence ``diagrams'' $\leftrightarrow$ ``irreducible representations''
also looks rather unnatural. This approach is traditional, and one can find it
in almost all textbooks and monographs on the subject,
for example, in one of the last books \cite{18}.
Under this approach,
considerable efforts are required to obtain any explicit
formula for characters of $S_n$.

Both these ways are important as well as indirect; they rest upon
deep and nontrivial auxiliary constructions. There is a natural question:
whether one can arrive at the main combinatorial objects of the
theory (diagrams, tableaux, etc.) in a more direct and natural fashion?

   We believe that the representation theory
 of the symmetric groups must satisfy the following three conditions:

\smallskip
{\advance\parindent by 3mm
\item{\rm(1)} The symmetric groups form a natural chain
($S_{n-1}$ is embedded
into $S_n$), and the representation theory of these groups
should be constructed inductively with respect to these embeddings,
that is, the representation theory  of $S_n$
should rely on the representation theory of $S_{n-1}$,
$n=1,2,\ldots$.
\item{\rm(2)} The combinatorics of Young diagrams and Young tableaux,
which reflects the branching rule for the restriction
$$
S_n\downarrow S_{n-1},
$$
should be introduced as a natural auxiliary element of the construction
rather than {\it ad hoc}; it should be deduced
from the intrinsic structure of the symmetric groups.
Only in this case the branching rule (which is one of the main theorems of
the theory) appears naturally  and not as
a final corollary of the whole theory.
\item{\rm(3)} The symmetric groups are Coxeter groups, and the methods of their
representation theory should apply to all classical series of Coxeter groups.

}
\smallskip
In this paper, we suggest a new approach, which satisfies the above
principles and makes the whole theory more natural and simple.
The following notions are very important for our approach:
\roster
\item Gelfand--Tsetlin algebra and Gelfand--Tsetlin basis
(GZ-algebra and GZ-basis);
\item Young--Jucys--Murphy (YJM-) elements;
\item algebras with a local system of generators (ALSG) as a general context
for the theory.
\endroster

The Gelfand--Tsetlin basis was defined by I.~M.~Gelfand and
M.~L.~Tsetlin in the fifties \cite{5, 6} for the
unitary and orthogonal groups. The general notion of GZ-algebra for inductive
limits of algebras can be introduced in the same way
 for an arbitrary inductive limit of semisimple algebras (this was done,
for example, in \cite{3}). For the general
definition of Gelfand--Tsetlin algebras and Young--Jucys--Murphy
generators, see also \cite{34}.

The notion of algebras or groups  with a local system
of generators and local relations (in short,
local algebras or groups)
generalizes Coxeter groups, braid groups, Hecke algebras,
locally free algebras, etc.\
(see \cite{30, 31}). This notion allows one to define an
inductive process of constructing representations,
which we apply here to the symmetric groups.

The special generators of the GZ-algebra of the symmetric group $S_n$
were essentially introduced  in  papers by A.~Young and then rediscovered
independently by A.-A.~A.~Jucys \cite{19} and G.~E.~Murphy \cite{24}.
These YJM-generators are as follows:
$$
\align
X_i &= (1\,i)+(2\,i)+\dots+(i-1\,i),\quad i=1,2,\dots,n;\\
X_0 &= 0,\quad X_1=(1,2), \quad \dots
\endalign
$$
There exist an invariant way to define them (see below),
which applies to a very general class of ALSG, in particular, to
all Coxeter groups. It is very important that these generators do not lie
in the centers of the corresponding group algebras, but nevertheless
generate the GZ-algebra, which contains all these  centers.

The complexity of the symmetric group
(compared, for example, to the full linear group) lies in the fact that
the Coxeter relations
$$
s_i s_{i+1} s_i = s_{i+1} s_i s_{i+1}
$$
for the generators $s_i$ of $S_k$ are not
commutation relations. Moreover, there is no sufficiently
large commutative subgroup of $S_k$ that
could play the role of a Cartan subgroup.
However, our approach in some way resembles
Cartan's highest weight theory, with the role of a Cartan subgroup played
by the  commutative GZ-subalgebra in $\C[S_n]$. The
Young--Jucys--Murphy generators of this subalgebra
diagonalize simultaneously
in any representation of $S_n$, and the whole
representation theory of $S_n$ is encoded in their spectrum.
The problem is, therefore, to describe this spectrum,
that is, to understand what eigenvalues of the YJM-elements can appear
and which of them appear in a given irreducible representation.

This problem is similar to the description of the dominant weights of a
reductive group.
We solve it using induction on $n$ and
elementary analysis of the commutation relation
$$
s_i X_i +1 = X_{i+1} s_i \,,\quad i=1,2,\dots,n-1,
\tag 0.1
$$
between the YJM-elements and the Coxeter generators $s_i$. In a sense,
the algebra $H(2)$ (the degenerate affine Hecke algebra
of order $2$) generated by $s_i$ and
two commuting elements $X_i$ and $X_{i+1}$ subject
to (0.1)
plays the same role in our paper as the group $\frak{gl}(2)$ plays
in the representation theory of
reductive groups.

Our exposition is organized as follows. We define the branching scheme of
irreducible representations of the symmetric groups $S_n$
and prove that it is a graph (rather than a multigraph), i.e., the multiplicities
of irreducible representations of $S_{n-1}$ in the restrictions of
irreducible representations of $S_n$ to $S_{n-1}$ are simple.
Then we study a maximal commutative subalgebra of the group algebra --- the
Gelfand--Tsetlin algebra, or the GZ-algebra, whose diagonalization in
each irreducible representation determines a linear basis of this
representation, and show that the spectrum of this algebra
is the set of integer vectors in $\Bbb R^n$ determined by simple
conditions described in Sec.~5 (so-called content vectors).
A vector satisfying these conditions is in turn just the vector
consisting of the
``contents'' of the boxes of a Young tableau (such a vector uniquely determines
the tableau), and thus we arrive at the main conclusion that
the bases of all irreducible complex representations of $S_n$
are indexed by Young tableaux. There is an equivalence relation on content
vectors: two vectors are equivalent if they belong to the same irreducible
representation. We prove that this equivalence of the
corresponding tableaux means that
they have the same Young diagram, and this completes the proof of the main
theorem --- the branching theorem: the branching graph (Bratteli diagram)
of irreducible representations of the symmetric groups $S_n$
coincides with the graph of Young diagrams (the Young graph).

Two facts allow us to realize this plan: first, we
choose the so-called
Young--Jucys--Murphy generators of the
Gelfand--Tsetlin algebra and consider the spectrum
with respect to these generators;
and, second, we can explicitly describe the representations of the degenerate
affine Hecke algebra $H(2)$, which plays the role of the ``increment''
in the inductive step from the group algebra $\C[S_{n-1}]$ to
the group algebra $\C[S_n]$. This step can be realized because of the role that
is played by the Coxeter generators of $S_n$ and the Coxeter
relations between them: they directly give conditions on elements of
the spectrum of the GZ-algebra (content vectors). One of the main
advantages of our construction of the representation theory of the symmetric
groups (and other series of Coxeter groups) is that {\it we obtain the
branching rule simultaneously with the description of representations, and
introduce Young diagrams and tableaux using only the analysis of the
spectrum of the GZ-algebra}. One may say that our plan also realizes a
noncommutative version of Fourier analysis on the symmetric groups,
in which
{\it the set of Young tableaux appears in a natural
way as the spectrum of a dual
object to $S_n$, and the set of diagrams gives the list of representations.}

As an application of these results, we derive the classical
Young formulas for the action of the Coxeter generators $s_i$ of $S_n$
and a new proof of the Murnaghan--Nakayama rule for the
characters of $S_n$. The final step in the proof
of the Young formulas is the same as in \cite{24};
in fact, the derivation of the Young formulas was Murphy's motivation for introducing
the elements $X_i$.  The novelty of our approach
compared to \cite{24} is that we
do not assume any knowledge of the representation theory
of $S_n$ and, on the contrary, construct the theory starting from
simple commutation relations.$^4$\footnotetext"$^4$"{From the viewpoint of the classical
representation theory of $S_n$, it may seem that using the whole
 inductive family,
$S_1\subset\dots\subset S_{n-1}\subset S_n$ to construct the representation
theory of the unique group $S_n$  is somewhat arbitrary (there are many
such families, although they are isomorphic). But it is this
``noninvariance'' that allows us to relate the theory to
Young diagrams and
tableaux; without it there is no branching theorem, no GZ-bases, no RSK
correspondence, etc. Moreover, without fixing
an inductive family, the correspondence ``irreducible
representations'' $\leftrightarrow$ ``Young diagrams''
 loses its precise sense and remains only an arbitrary act
of constructing the Specht modules. Of course, other inductive families (for
example, $S_2\subset S_4\subset\dots$ with periodic embeddings) lead to
other branching theorems and other bases.}

The first attempt to develop a new approach to the representation theory
of the symmetric groups was made in the papers \cite{30, 31},
where the notion of algebras with a local system of generators (ALSG)
was introduced.
The branching rule and Young's orthogonal form
were deduced in \cite{30}
from the Coxeter relation for the generators
of $S_n$ and the assumption that the branching
graph (see below) of $S_n$ is the Hasse diagram of a {\it distributive lattice}.
The approach presented in this paper does not require any additional
assumptions.

Our scheme can be applied to some other
ALSG, and first of all to the Coxeter groups of B--C--D series
and to wreath products of the symmetric groups with some finite groups.
All these generalizations will be considered elsewhere.

We do not attempt to give a complete bibliography on the subject.
Proper analogs of the Young--Jucys--Murphy elements for the
infinite symmetric group $S_\infty$ proved to be an
extremely powerful tool in infinite-dimensional representation
theory; see \cite{8, 9, 10, 11}. For the representation theory of
the infinite symmetric group, see
also \cite{20, 32, 3, 21}. In the series of papers \cite{30, 31, 32},
the first author develops a new approach to the representation theory of
$S_n$ in connection with asymptotic problems.

There are numerous other applications of the YJM-elements to classical
representation theory (see, for example, \cite{15}; we learned about this
important preprint after our paper was completed). The Young--Jucys--Murphy
elements arise naturally in connection with higher Capelli identities (see
\cite{27}). In \cite{13, 16},
these elements were considered in the context of the theory of degenerate
affine Hecke algebras. Young--Jucys--Murphy elements for Coxeter groups were
defined in \cite{26, 28}; among earlier papers, we mention \cite{7}.

In what follows, the reader is supposed to be familiar only with
elementary facts from the abstract representation theory of finite
groups. We will not use any facts from the representation theory of the
symmetric groups.

A short announcement of our results was made in \cite{4}.

\head
1. Gelfand-Tsetlin algebra and Gelfand--Tsetlin basis
\endhead

Consider an inductive chain of finite groups
$$
\{1\}=G(0)\subset G(1)\subset G(2) \subset \dots \,. \tag 1.1
$$
By $G(n)\w$ denote the set of
equivalence classes of
irreducible complex representations of the group $G(n)$.
By definition, the  {\it branching graph} (more precisely, the {\it branching multigraph}),
also called the {\it  Bratteli diagram},  of this chain is the following
directed graph.  Its vertices
 are the elements of the set (disjoint union)
$$
\bigcup_{n\ge 0} G(n)\w \,.
$$
Denote by $V^{\lambda}$ the $G(n)$-module corresponding to a representation
${\lambda}\in G(n)\w$.
Two vertices
$\mu\in G(n-1)\w$ and $\l\in G(n)\w$ are joined by $k$
directed edges (from $\mu$ to $\lambda$) if
$$
k=\dim\Hom_{G(n-1)}(V^\mu,V^\l)\,,
$$
that is, if $k$ is the multiplicity of $\mu$ in the restriction
of $\l$ to the group $G(n-1)$.
We call the set $G(n)\w$ the $n$th {\it level} of the
branching graph.
We write
$$
\mu\ar\l
$$
if $\mu$ and $\l$ are connected by an edge in the branching graph; and
$$
\mu\subset\l \,,
$$
where $\mu\in G(k)\w$, $\l\in G(n)\w$, and $k\le n$,
if the multiplicity of $\mu$ in the restriction of
$\l$ to $G(k)$ is nonzero.
In other words, $\mu\subset\l$ if there is a path
from $\mu$ to $\l$ in the branching graph.
Denote by $\net$ the unique
element of $G(0)\w$. The same definition
of the branching graph applies to
any chain
$$
M(0)\subset M(1)\subset M(2)\subset \dots
$$
of finite-dimensional semisimple algebras (see \cite{3} and references
therein). If the multiplicities
of all  restrictions are equal $0$ or $1$, then this
diagram is a graph (and not  multigraph); in this case one says that the
{\it multiplicities are simple} or the  {\it branching is simple}.
It is well known,
and we will prove this in the next section, that this is the case for
the   symmetric groups $G(n)=S_n$ (see also, e.g., \cite{18, 17}).
If the branching is simple,  the decomposition
$$
V^\l=\bigoplus_{\mu\in G(n-1)\w,\,\mu\ar\l\,} V^\mu
$$
into the sum of irreducible $G(n-1)$-modules is canonical.
By induction, we obtain a canonical decomposition of
the module $V^\l$ into irreducible $G(0)$-modules
(i.e., one-dimensional subspaces)
$$
V^\l=\bigoplus_T V_T
$$
indexed by  all possible chains
$$
T=\l_0\ar\l_1\ar\dots\ar\l_n, \tag 1.2
$$
where $\l_i\in G(i)\w$
 and $\l_n=\l$.
Such chains are {\it increasing paths}
from $\net$ to $\l$ in the branching
graph (or multigraph).

Choosing a unit (with respect to the $G(n)$-invariant
inner product $(\cdot\,,\,\cdot)$ in $V^\l$)
vector $v_T$ in each one-dimensional space $V_T$, we obtain a basis $\{v_T\}$ in the module $V^\lambda$, which
is called the {\it Gelfand--Tsetlin basis (GZ--basis)}. In \cite{5, 6},
such a basis was defined
for representations of $SO(n)$ and $U(n)$;
we use the same term in the general situation (see \cite{3}).
By the definition of $v_T$,
$$
\C[G(i)]\cdot v_T, \quad i=1,2,\dots,n,  \tag 1.3
$$
is the irreducible $G(i)$-module $V^{\l_i}$.
It is also clear that $v_T$ is the unique (up to
a scalar factor) vector with this property.

By $Z(n)$ denote the center of $\C[G(n)]$.
Let $GZ(n)\subset\C[G(n)]$ be the algebra generated
by the subalgebras
$$
Z(1),Z(2),\dots,Z(n)
$$
Íof $\C[G(n)]$. It is readily seen that the algebra $GZ(n)$
is commutative. It is called the
{\it Gelfand--Tsetlin subalgebra (GZ-algebra)} of the inductive
family of (group) algebras. Recall the
following fundamental isomorphism:
$$
\C[G(n)] =\bigoplus_{\l\in G(n)\w} \End(V^\l) \, \tag 1.4
$$
(the sum is  over all equivalence classes of irreducible complex
representations).

\proclaim{Proposition 1.1} The algebra $GZ(n)$ is the algebra of
 all operators diagonal in the Gelfand--Tsetlin basis.
In particular, it is a maximal commutative subalgebra
of $\C[G(n)]$.
\endproclaim

\demo{Proof}
Denote by $P_T\in GZ(n)$  the product
$$
P_{\l_1} P_{\l_2} \dots P_{\l}, \quad P_{\l_i}\in Z(i)\,,
$$
of the central idempotents corresponding to the representations $\l_1,\l_2,\dots,\l$, respectively.
Clearly, $P_T$ is a projection onto $V_T$. Hence
$GZ(n)$ contains the algebra of operators
diagonal in the basis $\{v_T\}$, which is
a maximal commutative subalgebra of $\C[G(n)]$. Since $GZ(n)$ is
commutative, the proposition follows. \hfill{$\square$}
\enddemo

\noindent{\bf Remark 1.2.}
Note that by the above proposition, any vector from the Gelfand--Tsetlin
basis in any irreducible representation of $G(n)$ is
uniquely (up to a scalar factor) determined by the eigenvalues of the
elements of $GZ(n)$ on this vector.

\medskip\noindent{\bf Remark 1.3.}
For an arbitrary inductive family of semisimple algebras, the GZ-subalgebra
is a maximal commutative subalgebra if and only if the branching graph
has no multiple edges.
\medskip

The following criterion of  simple branching uses the important
notion of centralizer. Let $M$ be a semisimple finite-dimensional
$\C$-algebra, and let $N$ be its subalgebra;
the {\it centralizer} $\ZMN$ of this pair is the subalgebra of all elements of
$M$ that commute with $N$.

\proclaim{Proposition 1.4}
The following two conditions are equivalent.

\itemitem{\rm(1)}
The restriction of any finite-dimensional irreducible complex
representation of the algebra $M$ to $N$ has simple
multiplicities.
\itemitem{\rm(2)}
The centralizer $\ZMN$ is commutative.
\endproclaim

\demo{Proof}
Let $V^{\mu}$ and $V^{\lambda}$ be the finite-dimensional spaces
of irreducible representations of the algebras
$N$ and $M$, respectively. Consider the
$M$-module $\Hom_N(V^{\mu},V^{\lambda})$. It is an irreducible
$\ZMN$-module; thus it is one-dimensional if the centralizer
is commutative.

Conversely, if there exists an irreducible representation of the centralizer
$\ZMN$ of dimension
more than one, then the multiplicity of the restriction of some representation of
$M$ to $N$ is also greater than one.
\hfill{$\square$}\enddemo

In the next section, we will apply this criterion to the group algebras of
the symmetric groups.

\head
2. Young--Jucys--Murphy elements
\endhead

From now on we consider the case
$$
G(n)=S_n.
$$

First let us prove that the spectrum of the restriction
of an irreducible representation of
$S_n$ to $S_{n-1}$ is simple
(i.e., there are no multiplicities).
The proof reproduces the idea of the classical I.~M.~Gelfand's criterion
saying when a pair of groups --- a Lie group and its subgroup --- is
what was later called a Gelfand pair (this means that the subalgebra
consisting of those
elements of the group algebra that are biinvariant with respect to the
subgroup is commutative). We present this beautiful proof (I was reminded of
it by E.~Vinberg), because it is of very general character and uses the
specific features of the symmetric group as little as possible.

Recall (see Proposition~1.4) that the spectrum of the
restriction of a representation of a group to a subgroup is simple if and only if
the {\it centralizer} of the group algebra of the
subgroup in the group algebra of the whole group {\it is commutative}.

\proclaim{Theorem 2.1}
{The centralizer $Z(n-1,1)\equiv Z(\Bbb C[S_n],\Bbb C[S_{n-1}])$ of the subalgebra
$\Bbb C[S_{n-1}]$ in $\Bbb C[S_n]$ is commutative.}
\endproclaim

We begin with the following assertion.

\proclaim{Lemma 2.2}
{Every element $g$ of the symmetric group $S_n$ is conjugate to the inverse
element $g^{-1}$, i.e., there is $h\in S_n$ such that
$g^{-1}=hgh^{-1}$; moreover, the element $h$
can be chosen in the subgroup $S_{n-1}$.}
\endproclaim

\demo{Proof}
Indeed, it is obvious that for every
$k$ (in particular, for $k=n-1$), every permutation from
$S_{n-1}$ is conjugate to its inverse. Now let
$g \in S_n$; take a permutation
$h \in S_{n-1}\subset S_n$ that conjugates in
$S_{n-1}$ the permutation $g'\in S_{n-1}$ induced on
$1, \dots, n-1$ by $g$ (i.e., $g'=p_ng$, where $p_n$ is the virtual
projection; see the definition in Sec.~7) and its inverse
${g'}^{-1}$; that is, we take $h\in S_{n-1}$ such that
${g'}^{-1}=hg'h^{-1}$. Then
$h$, regarded as an element of $S_n$ with fixed point
$n$, realizes the desired conjugation:
$g^{-1}=hgh^{-1}$. Moreover,  we can choose $h$ to be an element of second order
with fixed point $n$.\hfill\qed
\enddemo

Recall a simple but important fact from the theory of involutive algebras.

\proclaim{Lemma 2.3} {An algebra over $\Bbb C$ with involution $*$
is commutative if and only if all its elements are normal (commute
with the its conjugate). If any real element of algebra is self
conjugate then algebra is commutative.}
\endproclaim

\demo{Proof} {Let $B$ be a $*$-algebra over $\Bbb C$, i.e., an algebra with a linear
anti-automorphism of second order all elements of which are normal. Then for any two self-adjoint
operators $a=a^*$ and $b=b^*$ we have:

$$(a+ib)(a-ib)=(a-ib)(a+ib),$$

which means that $a$ and $b$ commute and consequently the algebra is commutative.

The converse is trivial. Suppose $\Bbb C$-algebra is complex hall
of the algebra over $\Bbb R$ and $a$ and $b$ are real and so
self-conjugate, then there product also real and conjugate,
therefore
$$ab=(ab)^*=b^*a^*=ba,$$ and since real part and whole algebra is commutative.}\hfill\qed
\enddemo

Let us continue the proof of Theorem 2.1. For involutive algebras over $\Bbb C$ the statement is
slightly different: an algebra over $\Bbb C$ is commutative if and only if the involution is the
complex conjugation with respect to every realization of the algebra as the complex hall of a real
algebra.

\demo{Proof of Theorem 2.1}
As we have seen, it suffices to check that every real element of the
centralizer
$Z(n-1,1)\subset {\Bbb C}[S_n]$ is self-conjugate.

Let $$f=\sum_i c_i g_i,\qquad c_i \in \Bbb R,$$ be an arbitrary real element of
$Z(n-1,1)$; the above expansion is unique, because
$\{g,\; g \in S_n\}$ is a basis in ${\Bbb
C}[S_n]$. Since $f$ commutes with every $h$ from $S_{n-1}$, it follows
from the uniqueness of the expansion
$f=\sum_i c_i g_i$ that is does not change
if we apply an inner automorphism
$f \to hfh^{-1}$;
as we have proved above, we can choose
$h=h_i$ with $h_ig_ih_i^{-1}=g_i^{-1}$; then the summand
$c_ig_i$ turns into $c_ig_i^{-1}$. Thus, along with every summand
$c_ig_i$, the decomposition also contains the summand
$c_ig_i^{-1}$, which means that $f$ is a fixed point of
the anti-automorphism,
or that $f^*=f$.\hfill\qed
\enddemo

The analysis of the whole proof leads to the following statement.

\proclaim{Theorem 2.4}
{Let $A$ be a finite-dimensional $*$-algebra over $\Bbb R$, and let $B$ be its
$*$-subalgebra; assume that in
$A$ there is a linear basis
$G=\{g_i\}$ closed under the involution
(i.e., $G^*=G$), and for every
$i$ there exists an orthogonal ($b^*=b^{-1}$)
element $b_i \in B$ such that $b_ig_i{b_i}^*=g_i^*$. Then
the centralizer of the subalgebra
$B$ in the algebra $A$ is commutative, and thus the spectrum of the
restriction of irreducible representations of the algebra
$A$ to the subalgebra $B$ is simple.

If $A$ and $B$ are the group algebras of a finite group
$G$ and its subgroup $H\subset G$, respectively, and the basis consists
of elements of $G$, then this condition reads as follows: for every
$g \in G$, there exist elements
$h\in H$ and $g'\in G$ such that
$h^{-1}g'h=g^{-1}$; if we can take  $g'=g$, then we obtain
the above condition.}
\endproclaim

This criterion in the above form can be applied in many situations.
We emphasize that the above proof of the simplicity of spectrum
for the symmetric groups does not use in any way the analysis of representations of
$S_n$; and the fact itself is the first step towards the spectral analysis of
the symmetric groups and is based only on elementary algebraic properties of
the group. Later we will see that the simplicity of spectrum also
easily follows
from another fact concerning centralizers.

We will need not only the fact that the centralizer $Z(n-1,1)$ is commutative,
but also a more detailed description of this centralizer as well as its
relation to the Gelfand--Tsetlin algebra. We will describe the centralizer and the structure
of the Gelfand--Tsetlin algebra with the help of a special basis.

For $i=1,2,\dots,n$, consider the
following elements $X_i\in\C[S_n]$:
$$
X_i = (1\,i)+(2\,i)+\dots+(i-1\,\,i)
$$
(in particular, $X_1=0$). We will call
them the {\it Young--Jucys--Murphy elements} (or {\it YJM-elements}).

It is clear that
$$
X_i = \text{sum of all transpositions in $S_i$} \text{ \, --
\, }
\text{sum of all transpositions in $S_{i-1}$}, \tag 2.1
$$
that is,  $X_i$ is the difference of an element ofÍ $Z(i)$ and
an element ofÍ $Z(i-1)$.
Therefore $X_i\in GZ(n)$ for all $i\le n$.
In particular, the Young--Jucys--Murphy elements commute.

Let $A, B,\dots, C$ be elements or subalgebras of some algebra
$M$; by $\langle A, B,\dots, C\rangle$ denote the subalgebra of
$M$ generated by $A, B,\dots, C$.

\proclaim{Theorem 2.5}
In the algebra $\CSn$, consider its center $Z(n)$ and the center $Z(n-1)$
of the subalgebra $\C[S_{n-1}]\hookrightarrow\CSn$. Then
$$
Z(n)\subset \langle Z(n-1), X_n\rangle.
$$
\endproclaim

\demo{Proof}
Recall that
$$
X_n = \sum_{i=i}^{n-1}(i,n) = \sum_{i\neq j;\, i,j=1}^n (i,j) -
\sum_{i\neq j;\, i,j=1}^{n-1} (i,j).
$$
The second summand lies in $Z(n-1)$, hence the first one lies in $\langle
Z(n-1), X_n\rangle$. We have
$$
  X_n^2 = \sum_{i,j=1}^n(i,n)(j,n) = \sum_{i\neq j;\, i,j=1}^{n-1}
  (i,j,n) + (n-1)\I.
$$
Therefore the element $\sum_{i\neq j; i,j=1}^{n-1}(i,j,n)$
lies in $\langle Z(n-1),X_n\rangle$.  Adding  the element
$$
 \sum_{i\neq j\neq k;\, i,j,k=1}^{n-1}(i,j,k)
$$
from $Z(n-1)$, we obtain the following element from $Z(n)$:
$$
 \sum_{i\neq j\neq k;\, i,j,k=1}^n(i,j,k).
$$
Thus we have proved that the indicator of the conjugacy class
of cycles of length $3$ in $S_n$ also lies in
$\langle Z(n-1),X_n\rangle$.

Apply induction and consider the general case
$$
\align
  X_n\cdot \sum_{i_1,\dots,i_{k-1}=1}^n(i_1,\dots,i_{k-1},n) &=
  \sum_{i\neq i_s,\, s=1,\dots,n-1}(i,n)(i_1,\dots,i_{k-1},n)
 \\
  &\quad+\sum_{i,i_1,\dots,i_{k-1}}(i,i_1,\dots,i_{k-1},n).
\endalign
$$
Taking the sum of  the first summand with the class
$$
  \sum_{i,j,i_1,\dots,i_{k-1}=1}^n (i,j)(i_1,\dots,i_{k-1}),
$$
which lies in $Z(n-1)$, we obtain the conjugacy class in
$S_n$ of the product of a cycle of length $2$ with a cycle of length
$k$, i.e., an element from
$Z(n)$. Hence the second summand, the class of cycles of length
$k+1$, also lies in $\langle Z(n-1), X_n\rangle$.
Again taking its sum with the element
$$
  \sum_{i,i_1,\dots,i_k}(i,i_1,\dots,i_{k-1},i_k)\in Z(n-1),
$$
we obtain the conjugacy class of cycles of length
$k+1$ in $S_n$.

Thus the classes of all one-cycle$^5$\footnotetext"$^5$"{By
one-cycle permutations, we mean permutations with
one nontrivial cycle.}
permutations in $S_n$ lie in $\langle Z(n-1), X_n\rangle$.
It remains to apply the classical theorem saying that the center of the group
algebra $\CSn$ is generated by multiplicative generators ---
the classes
of one-cycle permutations. This theorem reduces to the assertion that
the power sums $\sum_{i=1}^n x_i^r \equiv p_r$ form a multiplicative
basis in the ring of symmetric functions
(\cite{23, Chap.~1}). Thus
$$
  Z(n)\subset \langle Z(n-1), X_n \rangle. \eqno{\square}
$$

\enddemo

\proclaim{Corollary 2.6}
The Gelfand--Tsetlin algebra is generated by the Young--Jucys--Murphy
elements:
$$
  GZ(n) = \langle X_1, X_2, \dots, X_n \rangle.
$$
\endproclaim

\demo{Proof}
By definition,
$$
  GZ(n)=\langle Z(1),\dots, Z(n)\rangle.
$$
Clearly, $GZ(2)=\C[S_2]=\langle X_1=0, X_2 \rangle = \C$.

Assume that we have proved that
$$
  GZ(n-1) = \langle X_1, \dots, X_{n-1} \rangle.
$$
Then we must prove that
$$
  GZ(n) = \langle GZ(n-1), X_n \rangle.
$$
The inclusion
$$
  GZ(n) \supset \langle GZ(n-1), X_n \rangle
$$
is obvious, hence it suffices to check that
$$
  Z(n) \subset \langle GZ(n-1), X_n \rangle.
$$
But Theorem~2.5 implies
$$
  Z(n) \subset \langle Z(n-1), X_n \rangle
       \subset \langle GZ(n-1), X_n \rangle.\eqno{\square}
$$
\enddemo

\medskip\noindent{\bf Remark 2.7.}
Note that the YJM-elements do not lie in the corresponding centers:
$X_k\not\in Z(k)$, $k=1,\dots,n$. It might seem natural to search for
a basis of $GZ(n)$ consisting of elements of  the centers
$Z(1),\dots,Z(n)$. However, it is a
``noncentral'' basis that turns out to be useful.

\proclaim{Theorem 2.8}
The centralizer $Z(n-1,1)\equiv Z(\CSn,\C[S_{n-1}])$ of the subalgebra
$\C[S_{n-1}]$ in $\CSn$ is generated by the center $Z(n-1)$
of $\C[S_{n-1}]$ and the element $X_n$:
$$
  Z(n-1,1) = \langle Z(n-1), X_n \rangle.
$$
\endproclaim

\demo{Proof}
A linear basis in the centralizer
$Z(n-1,1)$ is the union of a linear basis in
$Z(n-1)$ and classes of the form
$$
  \sum (i_1^{(1)},\dots,i_{k_1-1}^{(1)},n) (i_1^{(2)},\dots,i_{k_2}^{(2)})
       \dots (i_1^{(3)},\dots,i_{k_3}^{(3)}),
$$
where the sum is taken over distinct indices
$i_s^l$ that run over all numbers from
$1$ to $n-1$. But taking the sum of
such classes with the classes
$$
  \sum (i_1^{(1)},\dots,i_{k_1}^{(1)}) (i_1^{(2)},\dots,i_{k_2}^{(2)})
       \dots (i_1^{(3)},\dots,i_{k_3}^{(3)})
$$
(the sum is  over all indices from $1$ to $n-1$) from $Z(n-1)$,
as in the proof of Theorem~2.5, we obtain all classes from
$Z(n)$. Hence a linear basis of $Z(n-1,1)$ can be obtained as a linear
combination of elements of the bases of
$Z(n-1)$
and $Z(n)$, i.e.,
$$
  Z(n-1,1)\subset \langle Z(n-1), Z(n) \rangle.
$$
And since $Z(n) \subset \langle Z(n-1), X_n \rangle$ (by Theorem~2.5),
the theorem follows.
\hfill{$\square$}\enddemo

\proclaim{Theorem 2.9}
The branching of the chain $\C[S_1]\subset\dots\subset\CSn$ is simple,
i.e., the multiplicities of the restrictions of irreducible representations
of  $\CSn$ to $\C[S_{n-1}]$ equal $0$ or $1$.
\endproclaim

\demo{Proof}
Since the centralizer $Z(n-1,1)$ is commutative (because
$Z(n-1,1)\subset \langle Z(n-1), X_n \rangle$), it suffices to apply
the simplicity criterion from Proposition~1.4.
\hfill{$\square$}\enddemo

\proclaim{Corollary 2.10}
The algebra $GZ(n)$ is a maximal commutative subalgebra of
$\CSn$. Thus in each irreducible representation of
$S_n$, the Gelfand--Tsetlin basis is determined up to scalar factors.
\endproclaim

This basis is called the {\it Young basis}.
A.~Young considered it in representations, but could not describe it
as a global basis, since this requires the notions of
GZ-algebra and YJM-elements,
which were not known then.

\def\Sp{\text{\rm Spec}\,}

The Young basis is a common eigenbasis of
the YJM-elements. Let $v$ be a vector of this basis
in some irreducible representation; denote by
$$
\al(v)=(a_1,\dots,a_n)\in\C^n
$$
the eigenvalues of $X_1,\dots,X_n$ on $v$.
Let us call the vector $\al(v)$ the {\it weight} of $v$.
Denote by
$$
\Sp(n) = \{ \al(v),\, v \text{ belongs to the Young basis}\}
$$
the  spectrum of the YJM-elements. By Theorem~2.5 and
Remark~1.2, a point $\al(v)\in\Sp(n)$ determines $v$
up to a scalar factor.
It follows that
$$
|\Sp(n)| = \sum_{\l\in S_n\w} \dim\l \,.
$$
In other words, the dimension of the Gelfand--Tsetlin algebra is equal
to the sum of the dimensions of all pairwise  nonequivalent irreducible
representations.

By the definition of the Young basis, the set $\Sp(n)$ is
in a natural bijection with the set of all paths \thetag{1.2}
in the branching graph. Denote this correspondence
by
$$
T\mapsto\al(T), \quad \al\mapsto T_\al \,.
$$
Denote by $v_\al$ the vector (unique up to a nonzero scalar factor)
of the Young basis corresponding to a weight $\al$.
There is a natural equivalence relation $\sim$ on $\Sp(n)$.
Write
$$
\al \sim \b, \quad \al,\b\in\Sp(n)\,,
$$
if $v_\al$ and $v_\b$ belong to the same irreducible $S_n$-module,
or, equivalently, if the paths $T_\al$ and $T_\b$ have the same end.
Clearly,
$$
|\Sp(n)/\sim| = |S_n\w|\,.
$$
Our plan is to
\roster
\item
describe the set $\Sp(n)$,
\item
describe the equivalence relation $\sim$,
\item
calculate the matrix elements in the Young basis,
\item
calculate the characters of irreducible representations.
\endroster

\head
3. The action of generators and the algebra $H(2)$
\endhead

The Coxeter generators
$$
s_i=(i\,\,i+1),\quad i=1,\dots,n-1,
$$
of the group $S_n$
commute except for neighbors. In \cite{30}, such generators were
called  {\it local}. Here
 ``locality'' is understood as in physics; it
means  that remote generators commute and hence do not
affect each other. The locality manifests itself in the following property of
the Young basis.

\proclaim{Proposition 3.1}
For any vector
$$
v_T, \quad T=\l_0\ar\dots\ar\l_n, \quad \l_i\in S_i\w,
$$
and any $k=1,\dots,n-1$, the vector
$$
s_k\cdot v_T
$$
is a linear combination of the vectors
$$
v_{T'}, \quad T'=\l'_0\ar\dots\ar\l'_n, \quad \l'_i\in S_i\w,
$$
such that
$$
\l'_i = \l_i, \quad i\ne k \,.
$$
In other words, the action of $s_k$ affects only the $k$th
level of the branching graph.
\endproclaim

\demo{Proof}
Let $i>k$.
Since $s_k\in S_i$  and the module
$$
\C[S_i]\cdot v_T
$$
is irreducible, we have
$$
\C[S_i]\, s_k\cdot v_T = \C[S_i]\cdot v_T = V^{\l_i}\,, \tag 3.1
$$
where $V^{\l_i}$ is the irreducible $S_i$-module indexed by
$\l_i\in S_i\w$.

Since $s_k$ commutes with $S_i$, \thetag{3.1} also holds
for all $i<k$. Now it follows from \thetag{1.3} that $s_k\cdot v_T$
is a linear combination of the  desired vectors.
\hfill{$\square$}
\enddemo

In the same way it is easy to show
that the coefficients of this linear combination depend only
on $\l_{k-1},\l_k,\l'_k$, $\l_{k+1}$ and the choice of the scalar
factors in vectors of the Young basis. That is,
the action of $s_k$ affects only the $k$th
level and depends only on  levels  $k-1,k$, and $k+1$
of the branching graph. More precise formulas are given in Sec.~4.

We can also easily  deduce the above proposition   from the obvious
relations
$$
s_i X_j = X_j s_i, \quad j\ne i,i+1 \,. \tag 3.2
$$
The elements $s_i$, $X_i$, and $X_{i+1}$ satisfy a more
interesting (and well-known) relation
$$
s_i X_i + 1 =  X_{i+1} s_i\,, \tag 3.3
$$
which can obviously be rewritten as
$$
s_i X_i s_i + s_i =  X_{i+1} \,.
$$

The action of the YJM-elements on the Young basis
is also local. It readily follows from \thetag{2.1}
that if
$$
T=\l_0\ar\dots\ar\l_n
$$
and
$$
\al(T)=(a_1,\dots,a_n),
$$
then $a_k$ is the difference of a function of $\l_k$
and a function of $\l_{k-1}$ for all $k$.

Denote by $H(2)$ the algebra generated by the elements $Y_1$, $Y_2$, and $s$
subject to the following relations:
$$
s^2=1,\quad Y_1Y_2=Y_2Y_1,\quad sY_1+1=Y_2s \,.
$$
The generator $Y_2$ can be excluded, because $Y_2=sY_1s+s$,
so that the algebra $H(2)$ is generated by $Y_1$ and $s$,
but technically it is more convenient to include $Y_2$
in the list of generators.

This algebra will play the central role in what follows.
It is the simplest example of the {\it degenerate affine
Hecke algebra} (see below). It follows directly
 from these relations
 that {\it irreducible finite-dimensional representations of this
algebra are either one-dimensional or two-dimensional}. Indeed,
since $Y_1$ and $Y_2$ commute, they have a common eigenbasis;
taking any vector $v$ of this eigenbasis and applying the involution $s$ to $v$,
we obtain an $H(2)$-invariant subspace of dimension at most $2$. The
importance of the algebra $H(2)$ is based on the following obvious yet
useful fact.

\proclaim{Proposition 3.2}
The algebra $\CSn$ is generated by the algebra $\C[S_{n-1}]$ and the algebra
$H(2)$ with generators $Y_1=X_{n-1}$, $Y_2=X_n$, $s=s_n$, where
$X_{n-1}$ and $X_n$ are the corresponding YJM-elements and
$s_n=(n-1,n)$ is a Coxeter generator.
\endproclaim

Of course, the algebra $\CSn$ is generated by the subalgebra
$\C[S_{n-1}]$ and one generator $s_n$, but it is taking into account
the superfluous generators
$X_{n-1}$ and $X_n$ that allows us to use induction: each step from
$n-1$ to $n$ reduces to the study of  representations of
$H(2)$.

Another important property of the Coxeter generators and
YJM-elements is that the relations between them
are stable under shifts of indices. In \cite{30}, such relations
were called {\it stationary}.

\medskip\noindent{\bf Remark 3.3.}
The {\it degenerate affine Hecke algebra}
$H(n)$ is generated by commuting variables
$Y_1, Y_2,\dots,Y_n$ and Coxeter involutions
$s_1,\dots,s_{n-1}$ with relations ({3.2}),
({3.3}) (see \cite{16, 13}). If we put $Y_1 = 0$, then
the quotient of $H(n)$ modulo the corresponding ideal
is canonically isomorphic to $\CSn$.

\head
4. Irreducible representations of $H(2)$
\endhead

As already mentioned in Sec.~3, all irreducible representations of $H(2)$
are at most two-dimensional and have a vector $v$
such that
$$
Y_1 v = a v,\quad Y_2 v = b v,\quad a,b\in\C\,.
$$
If the vectors  $v$ and $s v$ are linearly independent,
then the relation
$$
s Y_1+1=Y_2 s \tag 4.1
$$
implies that $Y_1$ and $Y_2$ act in the basis $v,s v$
as follows:
$$
Y_1=
\left(\matrix
a&-1\\
0&b
\endmatrix
\right),
\quad
Y_2=
\left(\matrix
b&1\\
0&a
\endmatrix
\right), \quad
s=\left(\matrix
0&1\\
1&0
\endmatrix
\right)\,.
$$
If $b\ne a\pm1$, then this representation contains one of the two
one-dimensional subrepresentations;
denote it by $\pi_{a,b}$. If $b=a\pm1$, then this representation
contains the unique one-dimensional subrepresentation
$$
Y_1\mapsto a, \quad Y_1\mapsto b, \quad s_1\mapsto\pm1,
$$
in which $v$ and $sv$ are proportional; and, conversely,
if $v$ and $sv$ are proportional, then
$$
sv=\pm v,
$$
and ({4.1}) implies
$$
b=a\pm1.
$$

Note that always $a\ne b$, since otherwise
the operators $\pi_{a,b}(Y_i)$ cannot be diagonalized and thus such
representations cannot
occur in the action on the Young basis. If $a\ne b$,
then the operators
$\pi_{a,b}$ can be diagonalized, for example, as
follows:
$$
Y_1=
\left(\matrix
a&0\\
0&b
\endmatrix
\right),
\quad
Y_2=
\left(\matrix
b&0\\
0&a
\endmatrix
\right),
\quad
s =
\left(\matrix
\frac1{b-a}&1-\frac1{(b-a)^2}\\
1&\frac1{a-b}
\endmatrix
\right)\,. \tag 4.2
$$

Let us formulate our results as a proposition which describes representations
in terms of transformations of weights (i.e., eigenvectors).

\proclaim{Proposition 4.1}
Let
$$
\al=(a_1,\dots,a_i,a_{i+1},\dots,a_n)\in\Sp(n) \,.
$$
Then $a_i\in\Bbb Z$ and
\roster
\item $a_i\ne a_{i+1}$ for all $i$;
\item if $a_{i+1}=a_i\pm 1$, then $s_i \cdot v_\al = \pm v_\al$;
\item if $a_{i+1}\ne a_i\pm1$, then
$$
\al'=s_i\cdot\al=(a_1,\dots,a_{i+1},a_{i},\dots,a_n)\in\Sp(n)
$$
and $\al'\sim\al$ (see Sec.~2 for the definition of the
equivalence relation $\sim$). Moreover,
$$
v_{\al'}=\left(
s_i - \frac1{a_{i+1}-a_i}
\right) v_{\al},
$$
and
the elements $s_i,X_i,X_{i+1}$ act in the basis $v_\al,v_{\al'}$ by
formulas \thetag{4.2} with $Y_1$ replaced by $X_i$ and $Y_2$
replaced by $X_{i+1}$.
\endroster
\endproclaim

\def\Co{\text{\rm Cont}\,}
Recall that the transpositions $s_i$ from claim (3) of Proposition~4.1
are Coxeter transpositions. In order to emphasize their role in
the context of this section (as operations on weights $\al$),
we call them {\it admissible} transpositions.
Admissible transpositions preserve the set $\Sp(n)$ and the set $\Co(n)$
defined in the next section.
The two cases of this proposition correspond to the cases of
chain and square from Sec.~7.

Note that if $a_{i+1}\ne a_i\pm1$, then in the basis
$$
\biggl\{ v_{\al},c_i(s_i-d_i{\I})v_\al\biggr\},
  $$
where $c_i=(a_{i+1}-a_i)^{-1}$, $d_i=(1-c_i^2)^{-1/2}$, the matrix of
the transposition $s_i$ is orthogonal:
  $$
  s_i = \left(\matrix
  1/r&\sqrt{1-1/r^2}\\
  \sqrt{1-1/r^2}&-1/r
  \endmatrix\right),
  $$
where $r=a_{i+1}-a_i$. In Young's papers, this difference was called
the {\it axial} distance; it is the difference of the contents
(see Sec.~5) of the corresponding boxes of Young tableaux.

\head
5. Main theorems
\endhead

In this section, we describe the set $\Sp(n)$ introduced in Sec.~2 and
the equivalence relation $\sim$.
Let us introduce the set
$\Co(n)$ of {\it content vectors} of length $n$.

\proclaim{Definition}
We say that $\al=(a_1,\dots,a_n)$ is a content vector,
$$
\al=(a_1,\dots,a_n)\in\Co(n),
$$
if $\al$ satisfies the following conditions:
\smallskip
{\advance\parindent by 3mm
\item{\rm(1)} $a_1=0$;
\item{\rm(2)} $\{a_q-1,a_q+1\}\cap\{a_1,\dots,a_{q-1}\}\ne\net$
for all $q>1$ (i.e., if $a_q>0$, then $a_i=a_q-1$ for some $i<q$;
and if $a_q<0$, then $a_i=a_q+1$ for some $i<q$);
\item{\rm(3)} if $a_p=a_q=a$ for some $p<q$, then
$$
\{a-1,a+1\}\subset\{a_{p+1},\dots,a_{q-1}\}
$$
(i.e., between two occurrences of $a$ in a content vector there should
also be occurrences of $a-1$ and $a+1$).
}
\endproclaim

It is clear that
$$
\Co(n)\subset\Z^n \,.
$$

\proclaim{Theorem 5.1}
$$
\Sp(n)\subset\Co(n)\,. \tag 5.1
$$
\endproclaim

We need the following lemma.

\proclaim{Lemma 5.2}
Let
$
\al=(a_1,\dots,a_n)
$
and $a_i=a_{i+2}=a_{i+1}-1$ for some $i$, i.e., $\al$
contains a fragment of the form $(a,a+1, a)$. Then
$$
\al\notin\Sp(n)\,.
$$
\endproclaim

\demo{Proof}
Let $\al\in\Sp(n)$. By claim (2) of Proposition~4.1,
$$
s_i v_\al = v_\al, \quad s_{i+1} v_\al = - v_\al \,,
$$
i.e., $s_i s_{i+1} s_i v_{\al} = -v_{\al}$, but $s_{i+1} s_i
s_{i+1} v_{\al} = v_{\al}$,
 contradicting the Coxeter relations
$$
s_i s_{i+1} s_i = s_{i+1} s_i s_{i+1} \,. \eqno{\square}
$$
\enddemo

\demo{Proof of Theorem~5.1}
Let $\al=(a_1,\dots,a_n)\in\Sp(n)$. Since $X_1=0$, we have
$a_1=0$.

Let us verify conditions \therosteritem{2} and
\therosteritem{3} by induction on $n$. The case $n=2$ is trivial.
Assume now that $\{a_n-1,a_n+1\}\cap\{a_1,\dots,a_{n-1}\}=\net$.
Then the transposition of $n-1$ and $n$ is admissible and
$$
(a_1,\dots,a_{n-2},a_n,a_{n-1})\in\Sp(n)\,.
$$
Hence $(a_1,\dots,a_{n-2},a_n)\in\Sp(n-1)$ and, clearly,
$$
\{a_n-1,a_n+1\}\cap\{a_1,\dots,a_{n-2}\}=\net,
$$
 contradicting the induction hypothesis. This proves
the necessity of \therosteritem{2}.

Assume that $a_p=a_n=a$ for some $p<n$, and let
$$
a-1 \notin \{a_{p+1},\dots,a_{n-1}\}\,.
$$
We may assume that $p$ is  the largest possible, that is,
the number $a$ does not occur between $a_p$ and $a_n$:
$$
a \notin \{a_{p+1},\dots,a_{n-1}\}\,.
$$
Then, by the induction hypothesis, the number $a+1$ occurs in
the set $\{a_{p+1},\dots,a_{n-1}\}$ at most once.
Indeed, if it occurred at least twice, then, by the induction hypothesis,
the number $a$ would also occur.
Thus we have two
possibilities: either
$$
(a_{p},\dots,a_{n})=(a,\ast,\dots,\ast,a)\,,
$$
or
$$
(a_{p},\dots,a_{n})=(a,\ast,\dots,\ast,a+1,\ast,\dots,\ast,a)\,,
$$
where $\ast,\dots,\ast$ stands for a sequence of
numbers different from $a-1$, $a$, $a+1$.

In the first case, applying $n-p-1$ admissible transpositions,
we obtain
$$
\al\sim\al'=(\dots,a,a,\dots)\,,
$$
which contradicts claim (1) of Proposition~4.1.

In the second case, the same argument yields
$$
\al\sim\al'=(\dots,a,a+1,a,\dots)\,,
$$
which contradicts Lemma~5.2.
\hfill{$\square$}
\enddemo

We will need another equivalence relation. Write
$$
\al\ssim\b, \quad \al,\b\in\C^n\,,
$$
if $\b$ is an admissible  permutation
(the product of admissible transpositions) of the entries of $\al$.
Now we are ready for the appearance of Young diagrams and tableaux.
Namely, we will see that vectors from $\Co(n)$ are the {\it content
vectors} of Young tableaux.

Recall some definitions.
Denote by $\Y$ the Young graph (see Fig.~1).
\midinsert
$$
%\epsfbox{c:/tom/tom.rus/tom307/ve1.eps}
\epsfbox{figure1.eps}
$$
\vskip-0.5cm
\botcaption {Fig.~1. {\rm The Young graph.}}
\endcaption
%\endinsert
%\midinsert
$$
%\epsfbox{c:/tom/tom.rus/tom307/ve2.eps}
\epsfbox{figure2.eps}
$$
\vskip-0.5cm
\botcaption {Fig.~2. {\rm Contents of boxes.}}
\endcaption
\endinsert

By definition,
the vertices of $\Y$ are Young diagrams,
and two vertices $\nu$ and $\eta$ are joined
by a directed edge if and only if $\nu\subset\eta$ and
$\eta/\nu$ is a single box. In this case we
write $\nu\ar\eta$.
Given a box
$\square\in\eta$, the number
$$
c(\square)=x\text{-coordinate of }\square -
y\text{-coordinate of }\square
$$
is called the {\it content} of $\square$ (see Fig.~2).

By $\T(\nu)$ denote the set of paths in $\Y$ from
$\net$ to $\nu$; such paths are called
{\it standard tableaux} or {\it Young tableaux}.
A convenient way to represent a
path $T\in\T(\nu)$,
$$
\net=\nu_0\ar\dots\ar\nu_n=\nu,
$$
is to write the numbers $1,\dots,n$ in the
boxes $\nu_1/\nu_0,\dots,\nu_n/\nu_{n-1}$ of
$\nu_n$, respectively. Put
$$
\T(n)=\bigcup_{|\nu|=n}\T(\nu)\,.
$$
The following proposition can easily be checked.

\proclaim{Proposition 5.3}
Let
$$
T=\nu_0\ar\dots\ar\nu_n\in\T(n)\,.
$$
The mapping
$$
T\mapsto\left(c(\nu_1/\nu_0),\dots,c(\nu_n/\nu_{n-1})\right)
$$
is a bijection of the set of tableaux $\T(n)$ and the set
of content vectors $\Co(n)$ defined at the beginning of this section.
We have $\al\ssim\b$, $\al,\b\in\Co(n)$, if and only if
the corresponding paths have the same end,
that is, if and only if they are tableaux with the same diagram.
\endproclaim

\demo{Proof} The content vector of any standard Young tableau
obviously satisfies conditions
{(1)}, {(2)}, and  {(3)} of the definition
of a content vector, and these conditions
uniquely determine the tableau as a sequence of boxes of
the Young diagram.
\hfill{$\square$}\enddemo

In terms of Young tableaux, admissible transpositions
are transpositions of numbers from different rows
and columns.

\proclaim{Lemma 5.4}
Any two Young tableaux $T_1,T_2\in \T(\nu)$ with
diagram $\nu$ can be obtained from each other
by a sequence of admissible transpositions. In other words,
if $\al,\b\in\Co(n)$ and $\al\ssim\b$, then
$\b$ can be obtained from $\al$ by admissible
transpositions.
\endproclaim

\demo{Proof}
Let us show that by admissible transpositions we can
transform any Young tableau $T\in\T(\nu)$, $\nu=(\nu_1,\dots,\nu_k)$,
to the following tableau with the same diagram (and horizontal monotone
numeration):
$$
%\epsfbox{c:/tom/tom.rus/tom307/ve3.eps}
\epsfbox{figure3.eps}
$$
\noindent corresponding to the content vector
$$
{\al}(T^\nu)=(0,1,2,\dots,\nu_1-1,-1,0,\dots,\nu_2-2,-2,-1,\dots)
$$
from $\Co(n)$.
To this end, consider the last box of the last row of $\nu$.
Let $i$ be the number written in this box of $T$. Transpose
$i$ and $i+1$, then $i+1$ and $i+2$, $\dots$, and,
finally, $n-1$ and $n$. Clearly,
all these transpositions are admissible,
and we obtain a tableau with the number $n$ written in the last
box of the last row. Now repeat the same procedure for
$n-1,n-2,\dots,2$. \hfill{$\square$}
\enddemo

\proclaim{Corollary 5.5}
If $\al\in\Sp(n)$ and $\al\ssim\b$, $\b\in\Co(n)$, then
$\b\in\Sp(n)$ and $\al\sim\b$.
\endproclaim

\noindent{\bf Remark 5.6.}
Our chain of transpositions from the proof of Lemma~5.4, which
connects $T$ and $T^\nu$,
is minimal possible in the following sense. Denote by $s$ the permutation
that maps $T$ to $T^\nu$, i.e., that associates with the number
written in a given box of $T$ the number written in the same box of $T^\nu$.
Let $\ell(s)$ be the
number of inversions in $s$, that is,
$$
\ell(s)=\#\{(i,j)\in\{1,\dots,n\}\,|\,i<j,\,s(i)>s(j)\}\,.
$$
It is well known that $s$ can be written as
the product of $\ell(s)$ transpositions $s_i$ and
cannot be written as a shorter product$^6$
\footnotetext"$^6$"{Simply because $\ell(s_i g)=\ell(g)\pm 1$
for all $i$ and $g\in S_n$.}.
It is easy to see that our chain contains precisely
$\ell(s)$ admissible transpositions. In other words,
$\Co(n)$ is a ``totally geodesic'' subset of $\Z^n$
for the action of $S_n$. That is,  along with any two vectors $\Co(n)$
contains chains of vectors that realize the minimal path between them in the
sense of the word metric with respect to the Coxeter generators.

\medskip
In the proof of Lemma~5.4 we used the fact that we can transform every tableau
with a given diagram into any other tableau
with the same diagram using only Coxeter transpositions; it is this fact that
guaranteed that vectors of the Young basis with the same diagram
lie in the same representation. Thus with each irreducible representation we
can associate the structure of a graph, whose vertices are vectors of the
Young basis and edges are labelled by  Coxeter generators and connect
pairs of vectors that can be transformed into each other by the corresponding
generator. These graphs generalize the Bruhat graph (order) on the group
$S_n$.

\medskip\noindent{\bf Remark 5.7.}
The first author
(see \cite{2}) introduced  the so-called adic transformations
on the spaces of
paths of graded graphs; in particular,
the {\it Young transformation (automorphism)} on the space of
infinite tableaux (i.e., paths in the Young graph). This transformation
sends a tableau to the next tableau in the lexicographic
order on the set of tableaux with a given diagram. Hence any two finite tableaux
with the same diagram lie on the same orbit of the Young automorphism.
The interval of the  orbit  that passes
 through tableaux with a given diagram
starts from the tableau shown in the above figure
(with horizontal monotone
numeration) and ends by the tableau with
vertical monotone numeration.
But,  of course, these orbits are not geodesic,
unlike the above-defined  chain of transformations,
 which constitutes only
a part of an orbit.

\def\N{\Bbb Z}
\medskip
Recall that the Young graph $\Y$ is an infinite $\N$-graded graph
of Young diagrams with obvious grading and set of edges. The graph
consisting of the first $n$ levels is denoted by $\Y_n$.

We proceed to the proof of the central theorem of the paper.

\proclaim{Theorem 5.8}
The Young graph $\Y$ is the branching graph of the symmetric groups;
the spectrum of the Gelfand--Tsetlin algebra
$GZ(n)$ is the space of paths in the finite graph
$\Y_n$, i.e., the space of Young tableaux with $n$ boxes;
we have $\Sp(n)=\Co(n)$, where $\Sp(n)$ is the spectrum of
$GZ(n)$ with respect to the YJM-ge\-ne\-ra\-tors
$X_1,\dots,X_n$ and
$\Co(n)$ is the set of content vectors; the corresponding
equivalence relations coincide:
$\sim\, =\, \ssim$.
\endproclaim

\demo{Proof} As we have seen, the set of classes
$\Co(n)/\ssim$ is the set of classes of tableaux with the same diagram. Hence
$$
 \#\big\{\Co(n)/\ssim\big\} = p(n),
$$
where $p(n)$ is the number of partitions of the number
$n$, i.e., the number of diagrams with $n$ boxes. By Corollary~5.5,
each equivalence class in
$\Co(n)/\ssim$ either does not contain elements of the set
$\Sp(n)$, or is a subset of some class in
$\Sp(n)/\sim$. But
$$
\#\big\{\Sp(n)/\sim\big\} = \#\big\{S_n\w\big\} = p(n),
$$
because the number of irreducible representations is equal to the
number of conjugacy classes, which is again the number of partitions of $n$
(as the number of cycle types of permutations). Therefore each class of
$\Co(n)/\ssim$ coincides with one of the classes of
$\Sp(n)/\sim$. In other words,
$$
\Sp(n)=\Co(n) \quad\text{and}\quad \sim\,=\,\ssim .
$$
Obviously, it follows that the graph
$\Y$ is the branching graph of the  symmetric groups.
\hfill{$\square$}\enddemo

Thus the main theorem is proved. But the above analysis gives
much more than the proof of the branching theorem; in subsequent sections
we will use it to obtain an explicit model of representations
(Young's orthogonal form) and sketch the derivation of the formula
for characters.

\head
6. Young formulas
\endhead

Up to now we have been considering vectors $v_T$ of the Young basis
up to scalar factors. In this section, we will specify
the choice of these factors.

Let us start with the tableau $T^\l$ defined in the
proof of Lemma~5.4 (see the figure).
Choose any nonzero vector $v_{T^\l}$ corresponding to
this tableau.

Now consider a tableau $T\in\T(\l)$
and put
$$
\ell(T)=\ell(s)\,,
$$
where $s$ is the permutation that maps $T^\l$ to $T$.
Recall that $P_T$ denotes the orthogonal projection
onto $V_T$ (see Sec.~1). Put
$$
v_T= P_T\cdot s \cdot v_{T^\l}. \tag 6.1
$$
By Lemma~5.4, the permutation $s$ can be represented
as the product of $\ell(T)$ admissible transpositions.
Therefore, by definition \thetag{6.1} and formulas \thetag{4.2},
$$
s \cdot v_{T^\l} = v_{T} +
\sum_{R\in\T(\l),\, \ell(R)<\ell(T)} \gamma_R\, v_R,
\tag 6.2
$$
where $\gamma_R$ are some rational numbers. In particular,
assume that $T'=s_i T$ and
$$
\ell(T') > \ell(T)\,.
$$
Let
$$
\al(T)=(a_1,\dots,a_n)\in\Co(n)
$$
be the sequence of contents of boxes in $T$.
Then \thetag{4.2}, \thetag{6.1}, and \thetag{6.2} imply
$$
s_i \cdot v_T= v_{T'} + \frac1{a_{i+1}-a_{i}} v_T \,. \tag 6.3
$$
And, again by \thetag{5.2},
$$
s_i \cdot v_{T'}=\left(1-\frac1{(a_{i+1}-a_i)^2}\right) v_T -
\frac1{a_{i+1}-a_{i}} v_{T'} \,. \tag 6.4
$$
This proves the following proposition.

\proclaim{Proposition 6.1}
There exists a basis $\{v_T\}$ of $V^\l$ in which the
Coxeter generators $s_i$ act according to  formulas {\rm(6.3), (6.4)}.
All irreducible representations of $S_n$ are defined over
the field $\Q$.
\endproclaim

Another way to prove this proposition is to verify
directly that these formulas define a representation
of $S_n$ (that is, to verify the Coxeter relations).

The basis used above yields Young's {\it seminormal form} of $V^\l$.
If we normalize all vectors $v_T$, we
obtain Young's {\it orthogonal form} of $V^\l$.
This form is defined over $\R$. Denote
the normalized vectors by the same symbols $v_T$.
Then $s_i$ acts in the
two-dimensional space spanned by $v_T$ and $v_{T'}$ by
an orthogonal matrix. Thus
$$
s_i =
\left(
\matrix
r^{-1} & \sqrt{1-r^{-2}} \\
\sqrt{1-r^{-2}} & -r^{-1}
\endmatrix
\right) \,, \tag 6.5
$$
where
$$
r=a_{i+1}-a_{i}\,.
$$
This number is usually called the
{\it axial distance} (see \cite{18} and also \cite{30}).
If we write the action of the Coxeter generators
$s_i$ in the basis of standard tableaux, it looks as follows:
\smallskip
{\advance\parindent by 3mm
\item{$\bullet$}
if $i$ and $i+1$ are in the same row, then
$s_i$ leaves the tableau $T$ unchanged;
\item{$\bullet$}
if $i$ and $i+1$ are in the same column, then
$s_i$ multiplies  $T$ by $-1$;
\item{$\bullet$}
if $i$ and $i+1$ are in distinct rows and columns,
then in the two-dimensional space spanned by this tableau and the tableau
(which is also standard) in which the elements
$i$ and $i+1$ are swapped,
$s_i$ acts according to  ({6.5}).

}

\proclaim{Proposition 6.2}
There exists an orthogonal basis $\{v_T\}$ of $V^\l$ in which the
generators $s_i$ act according to  formulas \thetag{6.5}.
\endproclaim

\noindent{\bf Remark 6.3.}
Since the weight $\al(T^\l)$ of the vector $v_{T^\l}$ is the maximal
weight in $V^\l$ with respect to the lexicographic
order, we may call $\al(T^\l)$ the {\it highest weight} of $V^\l$
and call the vector $v_{T^\l}$ the  {\it highest vector} of $V^\l$.

\head{7. Comments and corollaries}
\endhead

The previous sections contain the construction of the first part
of the representation theory of the symmetric groups:  the description of
irreducible representations, branching of representations, expressions for
the Coxeter generators in representations. In particular, we have revealed
the intrinsic connection between the combinatorics of Young diagrams and
tableaux and the Young graph on the one hand and the representation theory of the
symmetric groups on the other hand.

The further plan, which includes studying
the relation to symmetric functions
(characteristic map), formulas for characters, the theory of induced
representations, the Littlewood--Richardson rule, the relation to the
representation theory of $GL(n)$ and Hecke algebras and to asymptotic theory,
can also be realized with the help of the same ideas, the main of which is an inductive
approach to the series of symmetric groups.

From all these topics, in the next section we will only sketch the proof of the
Murnaghan--Nakayama rule, leaving the rest to another occasion.

In this section, we will give several simple corollaries from the results
obtained in Secs.~1--6. First of all, we will deduce corollaries from the
branching theorem, which claims that the branching of irreducible
representations of the groups $S_n$ is described by the Young graph.

\proclaim{Corollary 7.1}
The multiplicity of an irreducible representation
$\pi_\mu$ of
$S_n$ in a representation $\pi_\la$ of $S_{n+k}$ is equal to the number
of paths between the diagrams
$\la$ and $\mu$ ($\la\partn n+k$,
$\mu\partn n$); in particular, if $\mu\not\subset\la$, it is equal to
$0$, and in the general case it does not exceed
$k!$, this estimate being sharp.
\endproclaim

\demo{Proof}
Only the last claim needs to be proved. The number of tableaux in the skew
diagram $\la/\mu$ does not exceed the number of different ways to
add $k$ new boxes successively to the diagram $\mu$.
If these $k$ boxes can be added to different rows and columns, the number of
these ways equals $k!$.
\hfill{$\square$}\enddemo

In particular, if $k=2$, then we have only three different cases:
\smallskip
{\advance\parindent by 3mm
\item{(1)}
the multiplicity of $\mu$ in ${\lambda}$ is equal to $0$, and the vertices
$\mu$ and
${\lambda}$ are not connected in the branching graph;
\item{(2)}
the multiplicity is equal to $1$, and the interval connecting
$\mu$ and
${\lambda}$ is the chain
$$
\mu\text{---}\nu\text{---}{\lambda};
$$
\item{(3)}
the multiplicity is equal to $2$, and the interval between $\mu$ and ${\lambda}$
is the square
$$
\mu\,{}^{\hbox{$\diagup^{\raise 1pt \hbox{$\,\nu\,$}}\diagdown$}}
_{\hbox{$\diagdown_{\raise -3pt
\hbox{$\,\eta\,$}}\diagup$}}\,{\lambda}.
$$

}

\def\ar{\nearrow}
\def\Z{\Bbb Z}

In the case of the chain, the transposition
$s_{l+1}$ multiplies all vectors of the form
$$
v_T, \quad T=\dots\ar\mu\ar\nu\ar{\lambda}\ar\dots,
$$
by a scalar (which is equal to
$\pm1$ in view of the relation $s_{l+1}^2=1$).
The action of the permutation $s_{l+1}$ in the case of the square
was considered in the previous section.

Note that the Young graph is the so-called Hasse diagram of the
distributive lattice of finite ideals of the lattice
$\Z_+ + \Z_+$, hence intervals in the Young graph have a standard
description, and the generic interval is a Boolean algebra.
{\it A priori} this important fact is quite nonobvious, but eventually it
proved to be a corollary of the Coxeter relations.
Taking it as an assumption, one can also derive the branching theory
(see \cite{30, 31}).

The next important conclusion is an abstract description of the
Young--Jucys--Murphy generators based on the previous results.

Define a mapping
$$
\tilde p_n:S_{n}\to S_{n-1}
$$
by the following operation of deleting the last symbol:
$$
\tilde p_n((\dots,n,\dots)(\dots)\dots(\dots))=
((\dots,\not{\!n},\dots)(\dots)\dots(\dots)),
$$
where the parentheses contain the cycle decomposition of a permutation
$g\in S_n$ and $\tilde p_n$ leaves
all cycles except the first one, which contains $n$,
unchanged and delete $n$ from the first cycle.
\def\I{\Bbb I}%
\def\CSn{\Bbb C[S_n]}%
\def\C{\Bbb C}%
The mapping $\tilde p_n$ enjoys the following
obvious properties:
\smallskip
{\advance\parindent by 3mm
  \item{(1)}
$\tilde p_n(\I_n)=\I_{n-1}$, where $\I_k$ is the identity in $S_k$;
  \item{(2)}
$\tilde p_n|_{S_{n-1}} = \text{id}_{S_{n-1}}$\quad ($S_{n-1}\subset
S_n$);
  \item{(3)}
$\tilde p_n(g_1 h g_2) = g_1 \tilde p_n(h) g_2$,\quad $g_1,g_2\in
S_{n-1}$, $h\in S_n$.

}

\smallskip
Note that conditions {(1)} and {(2)} follow from {(3)}. Indeed,
{(3)} implies
$$
  \tilde p_n(g\I)=g \tilde p_n(\I) = \tilde p_n(\I g)=\tilde p_n(\I)g
$$
for all $g\in S_{n-1}$, whence $\tilde p_n(\I)=\I$. But then
$\tilde p_n(g)=g$ for $g\in S_{n-1}$.

By $p_n$ denote the extension of the mapping
$\tilde p_n$ by linearity to the group algebra
$\CSn$:
$$
p_n:\C[S_{n}]\to\C[S_{n-1}].
$$
Thus $p_n$ is a projection of the algebra $\CSn$ to the subalgebra
$\C[S_{n-1}]$. For $n=2,3$, such a projection is not unique, but for
$n\ge 4$, condition {(3)} uniquely determines an operation
$\tilde p_n:S_n\to S_{n-1}$. It is easy to see that the existence of
$\tilde p_n$ means the existence of an
$S_{n-1}$-biinvariant partition of
$S_n$ into $(n-1)!$ sets of $n
$ elements. Such a property of a pair of groups
$(G,H)$ is not satisfied often; however, there exists a generalization of
this construction to
semisimple algebras (in particular, to group algebras)
in the most general case.

\proclaim{Proposition 7.2}
$$
p_n^{-1}(\{c\I\}) \cap Z(n-1,1)  =  \{aX_n + b\I\},\quad
a,b,c\in\C.
$$
In other words, the inverse image of scalars intersects the centralizer
of $S_{n-1}$ in $\CSn$ by the two-dimensional subspace spanned by
the identity and the Young--Jucys--Murphy element
$X_n$. In particular,
$X_n$ is uniquely determined (up to scalar) as an element of the intersection
$$
p_n^{-1}(\{c\I\}) \cap Z(n-1,1)
$$
that is orthogonal to constants.
\endproclaim

\demo{Proof}
If $p_n(\sum_{g\in S_{n-1}}c_g g) = c\I$, then the element $A =
\sum_{g\in S_{n-1}}c_g g$ must be a linear combination
of the form
$A=\sum_{i=1}^n b_i(i,n)$. Such an element commutes with
$S_{n-1}$ if and only if
$$%\begin{gather*}
  b_1 = \dots=b_{n-1}=a,\quad b_n=b,\quad\text{i.e., \ } \
   A  = aX_n + b\I.\eqno{\square}
$$

\enddemo

The projection $p_n$ allows us to define the inverse spectrum (projective
limit) of the groups
$S_n$ regarded as $S_{n-1}$-bimodules:
$$
  \lim\limits_{\leftarrow}(S_n,\tilde p_n)= \goth S;
$$
the space $\goth S$ is no longer a group, but on this space there is a
left and right actions of the group
$S_{\infty}$ of finite permutations, because the projection
$\tilde p_n$ commutes with the left and right actions of
$S_{n-1}$ for all $n$. In \cite{20}, this object was called
the {\it space of virtual permutations}; it is studied in detail in
\cite{21}. There is a generalization of this construction to other inductive
families of groups and algebras.

\def\Zlk{Z(l,k)}

In conclusion we generalize the theorem on the centralizer
$\Zlk$ of
$\CSn$ in $\C[S_{n+k}]$.

\proclaim{Theorem 7.3 {\rm[10]}}
The centralizer
$$
  \Zlk \equiv \C[S_{n+k}]^{\CSn}
$$
is generated by the center $Z(n)$ of $\CSn\subset\C[S_{n+k}]$,
the group $S_n$ permuting the elements $n+1,\dots,n+k$, and
the YJM-elements $X_{n+1},\dots,X_{n+k}$.
\endproclaim

The main case $k=1$ is proved in Sec.~2. The general case can be proved by
the same method.

Note that this method of proof is different from and simpler
than that suggested in
\cite{4} and \cite{10, 11}; namely, it turns out to be useful to consider
first the subalgebra
$\langle Z(n),X_{n+1},\dots,X_{n+k}\rangle$, as in Sec.~2.

\medskip\noindent{\bf Remark 7.4.}
The formulas that describe the action of the symmetric group in representations
associated with skew diagrams (i.e., with diagrams
equal to the difference of two true Young
diagrams one of which contains the other)
are similar to the formulas from Sec.~6.

\medskip
\def\Zlk{Z(l,k)}
\def\Hom{\text{\rm Hom}}

Indeed, let ${\lambda}$  be a partition of
$l+k$ and $\mu$ be a partition of $l$ with
$\mu\subset{\lambda}$. By $V^{{\lambda}/\mu}$ denote the
$\Zlk$-module
$$
V^{{\lambda}/\mu}=\Hom_{S_l}(V^\mu,V^{\lambda}).
$$
It is clear that this module has an orthonormalized Young basis
indexed by all Young tableaux with the skew diagram
${\lambda}/\mu$
(which is similar to the basis of the representation associated with an ordinary
Young diagram). In this basis, the generators
$$
X_{l+i},\quad i=i,\dots,k,
$$
of the algebra $\Zlk$ act by multiplication by the content of
the $i$th box, and the Coxeter generators of the subgroup
$S_k\subset\Zlk$ act according to formulas ({6.5}).

We use Theorem~7.3 in the proof of the formula for characters in the next
section.

\def\g{\gamma}
\def\Vlm{V^{{\lambda}/\mu}}
\def\clm{\chi^{\lambda/\mu}}

\head
8. Characters of the symmetric groups
\endhead

In this section, we give a sketch of the proof of the Murnaghan--Nakayama rule
for characters of the symmetric groups.
In contrast to the previous sections, we do not recall definitions of some
well-known notions. The key role in the proof is played by
Proposition 8.3 based on Theorem~7.3.

Recall that a Young diagram $\g$ is called a {\it hook} if
$\g=(a+1,1^b)$ for some $a,b\in\Z_+$\,. The number
$b$ is called the {\it height} of the hook $\g$. Recall also
that a skew diagram $\l/\mu$ is called a
{\it skew hook} if it is connected and does not
have two boxes on the same diagonal. In other words, $\l/\mu$ is
a skew hook if the contents of all boxes of $\l/\mu$ form an
interval (of cardinality $|\l/\mu|$) in $\Z$. The number of
rows occupied by $\l/\mu$ minus 1 is called the
{\it height} of $\l/\mu$ and is denoted by
${\langle{\lambda}/\mu\rangle}$.
Put $k=|\l/\mu|$. Let $\Vlm$ be the representation of $S_k$
indexed by a skew diagram $\l/\mu$, and let
$\clm$ be the corresponding character. Our aim is to prove
the following well-known theorem.

\proclaim{Theorem 8.1}
There is the following formula:
$$
\chi^{\lambda/\mu}((12\ldots k))=
\cases
(-1)^{\langle\lambda/\mu\rangle}
 & \text{if \ } \lambda/\mu \text{ \ is a skew hook,}\\
  0 & \text{ otherwise}.
  \endcases\tag8.1
  $$

\endproclaim

Now suppose that $\rho$ is a partition of $k$.
Consider the following permutation from the
conjugacy class corresponding to
$\rho$:
$$
(12\dots\rho_1)(\rho_1+1\, \dots \,\rho_1+\rho_2)(\dots)\dots \,.
$$
It is clear that repeatedly applying the theorem to the action of
this permutation in the Young basis,
we obtain the following classical rule.

\def\lm{\lambda/\mu}

\proclaim{Murnaghan--Nakayama rule}
Let $\rho$ be a partition ofÍ $k$. The value
$\chi^{\lambda/\mu}_\rho$ of the character $\clm$ on a permutation
of cycle type $\rho$ equals
$$
\chi^{\lm}_\rho=\sum_S (-1)^{\raise 2pt\hbox{$\scriptstyle
\langle{S}\rangle$}},
$$
where the sum is over all sequences $S$,
$$
\mu=\l_0 \subset \l_1 \subset \l_2 \dots = \l\,,
$$
such that $\l_i/\l_{i-1}$ is a skew hook with $\rho_i$ boxes and
$$
{\langle S\rangle}= \sum_i \,\, {\langle
{\lambda}_i/{\lambda}_{i-1}\rangle}.
$$
\endproclaim

It is well known and can easily be proved (see, for
example, [23, Chap.~1, Ex.~3.11])
that
this rule is equivalent to all other descriptions
of the characters, such as
the relation
$$
p_\rho = \sum_\l \, \chi^\l_\rho \, s_\l \,
$$
for symmetric functions (see \cite{23})
or the determinantal formula (see \cite{23, 18}). Note
that the theorem we are going to prove is obviously a special case
of the Murnaghan--Nakayama rule.

The same proof of the
following proposition was also given in \cite{15}.

\proclaim{Proposition 8.2}
Formula \thetag{8.1} is true for $\mu=\varnothing$.
\endproclaim

\demo{Proof}
It is easy to see (for example, by induction; see also the
proof of Theorem~2.1)
that
$$
X_2 X_3 \dots X_k = \text{sum of all $k$-cycles in $S_k$}\,.
\tag 8.2
$$
The eigenvalue of \thetag{8.2} on any vector of the Young basis in $V^\l$
equals
$$
(-1)^b b!\, (k-b-1)!
$$
if $\l$ is a hook of height $b$, and vanishes otherwise.
Clearly, the number of one-cycle permutations in $S_k$ equals
$(k-1)!$, and
$$
\dim\l= \binom{k-1}b
$$
if $\l$ is a hook of height $b$. Taking the trace of \thetag{8.2}
in $V^\l$ proves the proposition. \hfill{$\square$}
\enddemo

\def\CSk{\Bbb C[S_k]}

\proclaim{Proposition 8.3}
For any vector $v$ from the Young basis of $\Vlm$,
$$
\CSk \cdot v = \Vlm.
$$
\endproclaim

\demo{Proof}
The space $\Vlm$ is an irreducible module over the degenerate affine
Hecke algebra $H(k)$.
The vector $v$ is a common eigenvector for all $X_i$.
Thus, by the commutation relations in $H(k)$, the space
$$
\CSk \cdot v
$$
is $H(k)$-invariant and hence equals $\Vlm$. \hfill{$\square$}
\enddemo

\proclaim{Proposition 8.4}
If $\lm$ is not connected, then
$$
\clm((12\dots k))=0 \,.
$$
\endproclaim

\def\Ind{\text{\rm Ind}}

\demo{Proof}
Assume that $\lm=\nu_1\cup\nu_2$, where $\nu_1$ and $     \nu_2$ are
two skew Young diagrams that have no common edge.
Let $a=|\nu_1|$, $b=|\nu_2|$.
Consider the subspace of $V^{\lm}$ spanned by all
tableaux of the form
$\l/\mu$ that have the numbers $1,2,\dots,a$ in the diagram $\nu_1$ and
the numbers  $a+1,\dots,k$ in the diagram $\nu_2$.
Obviously, the numbers of such tableaux equal precisely the number
of tableaux of the form $\nu_1$ and $\nu_2$, respectively.
Consider the action of the
subgroup $S_a\times S_b$
of $S_k$ on this subspace. It follows from the Young formulas
that it is isomorphic, as an $S_a\times S_b$-module, to
$$
V^{\nu_1}\otimes V^{\nu_2} \,.
$$
By Proposition 8.3, we have an epimorphism
$$
\Ind_{S_a\times S_b}^{S_k} V^{\nu_1}\otimes V^{\nu_2}
\longrightarrow \Vlm .\tag8.3
$$
The dimensions of both sides
of \thetag{8.3} equal
$$
\binom ka\, \dim\nu_1 \,\dim\nu_2 \,.
$$
Hence \thetag{8.3} is an isomorphism.

In the natural basis of the induced representation,
the matrix of the operator corresponding to the permutation $(12\dots k)$
(as well as any other permutation that is not conjugate
to any element of $S_a\times S_b$) has only zeros
on the diagonal. This proves the proposition. \hfill{$\square$}
\enddemo

\proclaim{Proposition 8.5} If $\lm$ has two
boxes on the same diagonal, then
$$
\clm((12\dots k))=0 \,.
$$
\endproclaim

\demo{Proof} Assume that there are two such boxes.
Then there is a diagram $\eta$ such that
$$
\mu\subset\eta\subset\l
$$
and $\eta/\mu$ is a $2\times2$ square
$$
\eta/\mu=\boxplus\,.
$$
That is, $\Vlm$ contains an $S_4$-submodule $V^\boxplus$.
By Proposition 8.3, we have an epimorphism
$$
\Ind_{S_4}^{S_k} V^{\boxplus}
\longrightarrow \Vlm \,. \tag 8.4
$$
By the branching rule and Frobenius reciprocity,
the left-hand side of \thetag{8.4} contains only
irreducible $S_k$-modules $V^\delta$ with $\boxplus\subset\delta$.
In particular, $\delta$ cannot be a hook, so that
$$
\chi^\delta((12\dots k))=0
$$
by Proposition 8.2.
This proves the proposition. \hfill{$\square$}
\enddemo

In fact, we have proved that under the assumptions of
Proposition 8.5,
$$
\Hom_{S_k}(V^\g,\Vlm) = 0
$$
for all hook diagrams $\g$.

\proclaim{Proposition 8.6} Assume that $\l/\mu$ is a skew hook.
Then for any hook $\g=(a+1,1^b)$,
$$
\Hom_{S_k}(V^\gamma,\Vlm) =\left\{
\aligned &\Bbb C,\quad
b = {\langle \lambda/\mu\rangle}, \\
&0 \quad\text{ otherwise}.
\endaligned\right.
$$
\endproclaim

\demo{Proof}
Since translations of a skew diagram obviously preserve
the corresponding $S_k$-module, we may assume that
$\l$ and $\mu$ are  minimal, that is,
$$
\l_1>\mu_1,\quad \l'_1>\mu'_1\,.
$$
Let us show that if $b<{\langle \lm\rangle}$, then
$$
\Hom_{S_k}(V^\g,\Vlm) = 0 \,.
$$
Indeed, the module $V^\g$ contains a nonzero $S_{k-b}$-invariant
vector, and $\Vlm$ contains no such vectors, because
there are no such vectors even in $V^\l$
(this follows from the branching rule).
The case $b>{\langle \lm\rangle}$ is similar.

Now assume that $b={\langle \lm\rangle}$. Consider the space
$$
\Hom_{S_k}(V^\g,V^\l) \,.
$$
It is easy to see, for example, from the following picture
(and the Young formulas)
$$
%\epsfbox{c:/tom/tom.rus/tom307/ve4.eps}
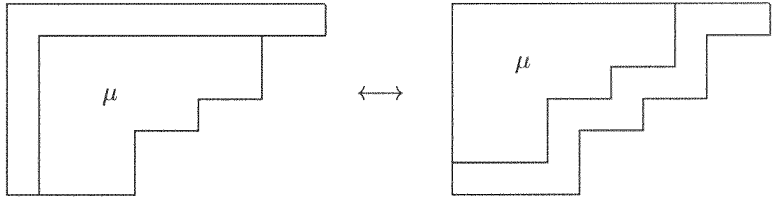
\epsfbox{figure4.eps}
$$
that this space is the irreducible $S_{|\mu|}$-module $V^\mu$.
Therefore
$$
\Hom_{S_k\times S_{|\mu|}}(V^\g\otimes V^\mu,V^\l)=\C \,,
$$
whence
$$
\Hom_{S_k}(V^\g,\Vlm) = \C.\eqno{\square}
$$
\enddemo

The theorem obviously follows from the propositions proved above.

\head{Acknowledgments}
\endhead

We would like to thank M.~Nazarov for useful information about the
literature, S.~Kerov and G.~Olshanski for helpful discussions, and the
referee for his valuable remarks.

The paper is supported by the  President of Russian Federation  grant for support of leading
scientific schools NSh-2251.2003.1.

\medskip
Translated by  the authors and N. Tsilevich.
\Refs

\item{1.}
 H.~Weyl,
{\it  The Classical Groups, Their Invariants and Representations},
Princeton Univ. Press, Princeton, New Jersey (1939).

\item{2.}
 A.~M.~Vershik,
``Uniform algebraic approximations of shift and multiplication
operators,''
{\it  Sov. Math. Dokl.},
{\bf  24}, 97--100 (1981).

\item{3.}
 S.~Kerov  and A.~Vershik,
``Locally semisimple algebras. Combinatorial theory and the
$K\sb0$-functor,''
{\it  J. Sov. Math.},
{\bf  38},
  1701--1733  (1987).

\item{4.}
 A.~M.~Vershik and A.~Yu.~Okounkov,
``An inductive method of presenting the theory of representations
of symmetric groups,''
  {\it  Russian Math. Surveys},
  {\bf   51},    No.~ 2,      355--356 (1996).
\item{5.}
 I.~M.~Gelfand and M.~L.~Tsetlin,
``Finite-dimensional representations of the group
of unimodular matrices,''
{\it  Dokl.\ Akad.\ Nauk SSSR},
{\bf   71},   825--828 (1950).
\item{6.}
 I.~M.~Gelfand and M.~L.~Tsetlin,
``Finite-dimensional representations of the group
of orthogonal matrices,''
{\it  Dokl.\ Akad.\ Nauk SSSR},
{\bf   71},   1017--1020  (1950).
\item{7.}
 V.~F.~Molchanov,
``On matrix elements of irreducible representations
of the symmetric group,''
{\it  Vestnik Mosk.\ Univ.}, {\bf   1},
    52--57 (1966).
\item{8.}
 A.~Yu.~Okounkov,
``Thoma's theorem and representations of the infinite bisymmetric group,''
{\it  Funct.\ Anal.\ Appl.},  {\bf   28}, No.~2, 101--170 (1994).

  \item{9.}
   A.~Yu.~Okounkov,
  ``On representations of the infinite symmetric group,''
  {\it  Zap. Nauchn. Semin. POMI}, {\bf 240}, 166--228 (1996).

  \item{10.}
   G.~I.~Olshanski,
``Extension of the algebra $U(\mathfrak g)$ for infinite-dimensional
classical Lie algebras $\mathfrak g$ and the Yangians $Y(\mathfrak{gl}(m))$,''
  {\it  Soviet Math. Dokl.}, {\bf  36},  569--573 (1988).

  \item{11.}
   G.~I.~Olshanski,
``Unitary representations of $(G,K)$-pairs that are connected
  with the infinite symmetric group $S(\infty)$,''
    {\it  Leningrad Math. J.}, {\bf   1},  No.~4, 983--1014 (1989).

\item{12.}
 I.~A.~Pushkarev,
``To the representation theory of wreath products of
finite groups with symmetric groups,''
  {\it  Zap. Nauchn. Semin. POMI}, {\bf 240}, 229--244 (1996).
\item{13.}
 I.~Cherednik,
``On special bases of irreducible
finite-dimensional representations of the
degenerate affine Hecke algebra,''
{\it  Funct.\ Anal.\ Appl.}, {\bf   20},   No.~ 1,
   87--88 (1986).
\item{14.}
 I.~Cherednik,
``A unification of Knizhnik--Zamolodchikov and Dunkl operators via
  affine Hecke algebras,''
  {\it  Invent. Math.}, {\bf   106},  No.~2, 411--431 (1991).

\item{15.}
 P.~Diaconis and C.~Greene,
``Applications of Murphy's elements,''
Stanford University Tech. Report, No.~335 (1989).
\item{16.}
 V.~Drinfeld,
``Degenerated affine Hecke algebras and Yangians,''
{\it  Funct. Anal. Appl.}, {\bf   20},  56--58 (1986).
\item{17.}
 G.~D.~James,
{\it  The Representation Theory of the Symmetric Group},
 Springer-Verlag,  Berlin (1978).
\item{18.}
 G.~James and A.~Kerber,
{\it  The Representation Theory of the Symmetric Group},
  Addison-Wesley,  Reading (1981).

\item{19.}
 A.~Jucys,
``Symmetric polynomials and the center of the symmetric group ring,''
{\it  Reports Math. Phys.}, {\bf   5},  107--112 (1974).
\item{20.}
 S.~Kerov, G.~Olshanski, and A.~Vershik,
``Harmonic analysis on the infinite symmetric group. A deformation of
  the regular representation,''
{\it  C. R. Acad. Sci. Paris Ser. I. Math.}, {\bf   316},
 No.~8,  773--778 (1993).
\item{21.}
 S.~Kerov, G.~Olshanski, and A.~Vershik,
``Harmonic analysis on the infinite symmetric group,''
{\it   Invent. Math.},  to appear.

\item{22.}
 G.~Lusztig,
``Affine Hecke algebras and their graded version,''
{\it  J. Amer. Math. Soc.}, {\bf   2},  No.~3, 599--635 (1989).
\item{23.}
 I.~G.~Macdonald,
{\it  Symmetric Functions and Hall Polynomials}, 2nd edition,
 Oxford University Press,  New York (1995).
\item{24.}
 G.~Murphy,
``A new construction of Young's seminormal representation of the
  symmetric group,''
{\it  J. Algebra}, {\bf   69},  287--291 (1981).
\item{25.}
 M.~L.~Nazarov,
``Projective representations of the infinite
  symmetric group,''
 in: {\it Representation Theory and Dynamical Systems},
 A.~Vershik (ed.),
  Amer. Math. Soc., Providence, Rhode Island (1992).
\item{26.}
 M.~Nazarov,
``Young's orthogonal form for Brauer's centralizer algebra,''
{\it  J. Algebra}, {\bf   182},  No.~3, 664--693 (1996).

\item{27.}
 A.~Okounkov
``Young basis, Wick formula, and higher Capelli identities,''
{\it  Internat. Math. Res. Notices}, {\bf   17},  817--839 (1996).
\item{28.}
 A.~Ram,
``Seminormal representations of Weyl groups and Iwahori--Hecke algebras,''
{\it  Proc. London Math. Soc.} (3), {\bf   75},  99--133 (1997).

\item{29.}
 J.~D.~Rogawski,
``On modules over the Hecke algebra of a $p$-adic group,''
{\it  Invent. Math.}, {\bf  79},  No.~3, 443--465 (1985).
\item{30.}
 A.~Vershik,
``Local algebras and a new version of Young's orthogonal form,''
in: {\it  Topics in Algebra, Part $2$:
 Commutative Rings and Algebraic Groups},
 Banach Cent. Publ., Vol. 26, part 2 (1990),
pp.~  467--473.
\item{31.}
 A.~Vershik,
``Local stationary algebras,''
{\it   Amer. Math. Soc. Transl.}  (2), {\bf   148},  1--13 (1991).

\item{32.}
 A.~Vershik,
``Asymptotic aspects of the representation theory of symmetric groups,''
{\it  Selecta Math. Sov.}, {\bf   11},  No.~2,  159--180 (1992).

\item{33.}
{\it  Asymptotic Combinatorics with Applications to
Mathematical Physics},  A.~M.~Vershik (ed.),
 Springer-Verlag,  Berlin--Heidelberg--New York
(2003).

\item{34.}
 A.~M.~Vershik,
``Gelfand--Tsetlin algebras, expectations,
inverse limits, Fourier analysis,''
 in: {\it Unity of Mathematics. Conference dedicated to I.~M.~Gelfand},
 to appear.

\item{35.}
 D.~Zagier,
   Appendix to the book S.~K.~Lando, A.~K.~Zvonkin,
Graphs on Surfaces and Their Applications,
 Springer,  Berlin (2004).

\endRefs

\enddocument